\documentclass[12pt,leqno]{article}

\usepackage{latexsym,amsmath,amssymb,amsthm,amsfonts}

\renewcommand{\epsilon}{\varepsilon}

\newcommand{\R}{\mathbb{R}}

\begin{document}

\author{Andriy V.\ Bondarenko}
\title{On a spherical code in the space of spherical harmonics}
\maketitle
\begin{abstract}
In this short note we propose a new method  for construction new nice arrangement
on the sphere $S^d$ using the spaces of spherical harmonic.
\end{abstract}
\section {Introduction}
This paper is inspired by classical book J. H. Conway and N. J. A. Sloane~\cite{CS} and recent paper of H. Cohn and A. Kumar \cite{CK}. The exceptional arrangement of points on the spheres are discussed there. Especially interesting are constructions coming from well known $E_8$ lattice and Leech lattice $\Lambda_{24}$. The main idea of the paper is to use these arrangements for construction new good arrangements in the spaces of spherical harmonics $\mathcal{H}_k^d$. Recently we have use dramatically the calculations in these spaces to obtain new asymptotic existence bounds for spherical designs, see~\cite{BV}. Below we need a few facts on spherical harmonics.
Let $\Delta$ be the Laplace operator in $\mathbb{R}^{d+1}$ $$
\Delta =\sum_{j=1}^{d+1}\frac{\partial^2}{\partial x_j^2}. $$ We
say that a polynomial $P$ in $\mathbb{R}^{d+1}$ is harmonic if
$\Delta P=0$. For integer $k\ge 1$, the restriction to $S^d$ of a
homogeneous harmonic polynomial of degree $k$ is called a
spherical harmonic of degree $k$. The vector space of all
spherical harmonics of degree $k$ will be denoted by
$\mathcal{H}_k^d$ (see~\cite{MNW} for  details). The dimension of
$\mathcal{H}_k^d$ is given by $$ \dim\,\mathcal{H}_k^d=
\frac{2k+d-1}{k+d-1}\binom{d+k-1}{k}.$$\\
Consider usual inner product in $\mathcal{H}_k^d$  $$\langle P,
Q\rangle:=\int_{S^d}P(x)Q(x)d\mu_d(x),$$ where $\mu_d(x)$ is normalized Lebesgue measure on the unit sphere $S^d$.
Now, for each point $x\in S^d$ there exists a unique
polynomial $P_x\in \mathcal{H}_k^d$ such that $$\langle
P_x,Q\rangle=Q(x) \;\;\mbox{for all}\;\;Q\in\mathcal{H}_k^d.$$
It is well known that $P_x(y)=g((x,y))$, where $g$ is a corresponding Gegenbauer polynomial.
Let $G_x$ be normalized polynomial $P_x$, that is $G_x=P_x/g(1)^{1/2}$. Note that
$\langle G_{x_1}, G_{x_2}\rangle=g((x_1,x_2))/g(1).$  So, if we have some arrangement $X=\{x_1,\ldots,x_{N}\}$ on $S^d$ with known distribution of inner products $(x_i,x_j)$, then for each $k$ we have corresponding set $G_X=\{G_{x_1},\ldots,G_{x_{N}}\}$ in $\mathcal{H}_k^d$, also with known distribution of inner products. Using this construction we will obtain in the next section the optimal antipodal spherical $(35, 240, 1/7)$ code from minimal
vectors of $E_8$ lattice. Here is a definition.

{\bf Definition 1.} The antipodal set $X=\{x_1,\ldots,x_{N}\}$ on $S^d$
is called antipodal spherical $(d+1, N ,a)$ code, if $|(x_i,x_j)|\le a$, for some $a>0$ and for
all $x_i, x_j\in X$, $i\ne j$, which are not antipodal. Such code is called optimal if for any
antipodal set $Y=\{y_1,\ldots,y_{N}\}$ on $S^d$ there exists $y_i, y_j\in Y$, $i\ne j$, which are not antipodal
and $|(y_i, y_j)|\ge a$. \\ In the other words, antipodal spherical $(d+1, N ,a)$ code is optimal if $a$ is a minimal
possible number for fixed $N,d$.

\section {Construction and the proof of optimality}
Let $X=\{x_1,\ldots,x_{120}\}$ be any subset of $240$ normalized minimal vectors of $E_8$ lattice,
such that no pair of antipodal vectors presents in $X$. Take in the space $\mathcal{H}^7_2$ the polynomials
$$
G_{x_i}(y)=g_2((x_i,y)),  \qquad i=\overline{1,\ldots,120},
$$
where $g_2(t)=\frac 87t^2-\frac 17$ is a corresponding normalized Gegenbauer polynomial.
Since $(x_i,x_j)=0$ or $\pm 1/2$, for $i\ne j$, then $\langle G_{x_i}, G_{x_j}\rangle=g_2((x_i,x_j))=\pm 1/7$!
It looks really like a mystery the fact that $|g_2((x_i,x_j))|=const$, for any different $x_i, x_j\in X$. But exactly this is essential for the proof of optimality of our code. Since, $dim \mathcal{H}^7_2=35$, then the points
$G_{x_1},\ldots,G_{x_{120}},-G_{x_1},\ldots,-G_{x_{120}}$ provide antipodal spherical $(35, 240, 1/7)$ code.
Here is a proof of optimality. Take arbitrary antipodal set of points $Y=\{y_1,\ldots,y_{240}\}$ in $\R^{35}$. Then, the inequality
$$
\frac 1{240^2}\sum_{i,j=1}^{240}(y_i,y_j)^2\ge 1/35,
$$
implies that $(y_i,y_j)^2\ge1/49$, for some $y_i,y_j\in Y$, $i\ne j$, which are not antipodal. This immediately gives us an optimality of our construction. The other reason why it works, that is our set is also spherical $3$-design in $\R^{35}$.  We are still not able generalize this construction even for Leech lattice $\Lambda_{24}$. We also don't know whether the construction described above is an optimal spherical $(35, 240, 1/7)$ code.\\
{\bf Acknowledgement.} The author would like to thank Professor
Henry Cohn for the fruitful discussions on the paper.

\end{document}